\documentclass[11pt]{article}

\usepackage{amsfonts,amsmath,amssymb,cite}
\usepackage{amsmath,amsthm}
\usepackage{enumerate}
\usepackage{latexsym}
\usepackage{color}
\usepackage{tikz}
\usetikzlibrary{decorations.pathmorphing}
\usetikzlibrary{shapes.geometric}
\usetikzlibrary{svg.path}

\newtheorem{thm}{Theorem}
\newtheorem{ob}{Observation}
\newtheorem{lem}[thm]{Lemma}

\newtheorem{prop}[thm]{Proposition}
\newtheorem{ques}{Question}
\newtheorem{prob}{Problem}
\newtheorem{claim}{Claim}

\newcommand{\ong}{N_{\rm o}}
\newcommand{\opack}{\rho^{\rm o}}
\newcommand{\Lopack}{\rho^{\rm o}_L}

\newcommand{\cU}{{\cal U}}
\newcommand{\cF}{{\cal F}}

\begin{document}

\title{On graphs having one size of maximal open packings}
\author{
Bert L.~Hartnell$^{a}$ \and Douglas F.~Rall $^{b}$}

\maketitle

\begin{center}
$^a$ Saint Mary's University, Halifax, Nova Scotia, Canada\\
\medskip

$^b$ Department of Mathematics, Furman University, Greenville, SC, USA\\
\end{center}

\begin{abstract}
A set $P$ of vertices in a graph $G$ is an open packing if no two distinct vertices in $P$ have a common
neighbor.  Among all maximal open packings in $G$, the smallest cardinality is denoted $\Lopack(G)$ and
the largest cardinality is $\opack(G)$.  There exist graphs for which these two invariants are arbitrarily
far apart.  In this paper we begin the investigation of the class of graphs that have one size of maximal open packings.
By presenting a method of constructing such graphs we show that every graph is the induced subgraph of a graph in this class.
The main result of the paper is a structural characterization of those $G$ that do not have a cycle of order less than
$15$ and for which $\Lopack(G)=\opack(G)$.

\end{abstract}

{\small \textbf{Keywords:} open packing; well-covered} \\
\indent {\small \textbf{AMS subject classification:} 05C69, 05C75}
\section{Introduction} \label{sec:intro}

A subset $P$ of the vertex set of a graph $G$ is an \emph{open packing} in  $G$ if the open neighborhoods of vertices in $P$ are pairwise disjoint.  That is,
no pair of distinct vertices in $P$ have a common neighbor.  By a \emph{maximal open packing} we mean an open packing that is maximal with respect to set containment.  The cardinality of a largest open packing in $G$ is called the \emph{open packing number} of $G$ and is denoted by $\opack(G)$.    The \emph{lower open packing number} of $G$, denoted $\Lopack(G)$, is the minimum cardinality of a maximal open packing in $G$.  To see that these two numbers can differ by an arbitrary amount consider the tree $T_n$ of order $4n+2$ obtained from the disjoint union of two stars $K_{1,n}$ by subdividing each edge once and then adding an edge to make their centers adjacent.  For this tree, $\Lopack(T_n)=2$ while $\opack(T_n)=2n+2$.

The study of maximal open packings in graphs was initiated by Henning and Slater~\cite{he-sl-1999}.  They determined the lower and upper open packing numbers for paths
and cycles, and in a series of results they established bounds for $\opack(G)$ and $\Lopack(G)$.  In particular, they proved that for a connected graph $G$ of order $n$ with maximum degree $\Delta$ and minimum degree $\delta$, we have $\frac{n}{\Delta(\Delta-1)+1} \leq \Lopack(G) \leq \opack(G) \le \frac{n}{\delta}$.  Bre\v{s}ar, Kuenzel, and Rall~\cite{br-ku-ra-2020+} investigated graphs with a unique maximum open packing and showed that recognition of this class of graphs is polynomially equivalent to the recognition of the graphs with a unique maximum independent set.  Furthermore, they gave a structural characterization of the class of trees $T$ with a unique open packing of
cardinality $\opack(T)$.  A number of other
researchers have studied open packings.  For example, see the following~\cite{he-1998b, he-1998, ha-sa-2015, ha-sa-2016, sa-sa-2017, mo-ma-2019, mo-sa-2019}.

Open packings are related to total domination since every open neighborhood has a nonempty intersection with any total dominating set.  (A set $S$ is a total dominating set if every vertex in $G$ is adjacent to at least one vertex in $S$.)  This implies that the size of a smallest total dominating set in any graph is at least as large as its
open packing number.  In the class of trees these invariants are equal.  In 2005, Rall~\cite{ra-2005} proved that if $T$ is a nontrivial tree, then the cardinality of a smallest total dominating set is $\opack(T)$.  This result is ``parallel'' to the theorem of Meir and Moon~\cite{me-mo-1975}  that showed the domination number of a tree is equal to the cardinality of a largest (closed) packing, which is a set of vertices whose closed neighborhoods are pairwise disjoint.

Henning and Slater~\cite{he-sl-1999} also considered the complexity of computing open packings.  They proved that the following decision problem
is NP-complete even when restricted to bipartite or to chordal graphs.

\begin{center}
	\fbox{\parbox{0.99\linewidth}{\noindent
			{\sc Open Packing}\\[.8ex]
			\begin{tabular*}{0.95\textwidth}{rl}
				{\em Instance:} & A graph $G$ and a positive integer $k$.  \\
				{\em Question:} & Does $G$ have an open packing of cardinality $k$? \\
			\end{tabular*}
	}}
\end{center}

Hamid and Saravanakumar~\cite{ha-sa-2015} posed as an open problem the characterization of those graphs $G$ such that $\Lopack(G)=\opack(G)$. The main goal
of this paper is to give a structural characterization of the subclass of this class of graphs that have no cycles of length less than $15$.  For convenience we denote
by $\cU$ the set of all graphs $G$ such that $\Lopack(G)=\opack(G)$.  Equivalently, $G \in \cU$ if and only if every maximal open packing in $G$ has cardinality $\opack(G)$.
Note that {\sc Open Packing} is solvable in linear time for the class $\cU$ since a greedy algorithm will always produce
an open packing of cardinality $\opack(G)$ for every $G \in \cU$.

Let $\cF$ be the family of all finite simple graphs, $G$, such that there is a weak partition $L, S_1,S_2,D_{11}, D_{12},D_2$
of $V(G)$ that satisfies the following properties.

\begin{enumerate}
\item $L$ is the set of leaves in $G$, and $S_1 \cup S_2$ is the set of support vertices in $G$.
\item The sets $S_1$ and $D_{11}$ are  independent, and $S_2$ induces a matching in $G$.
\item $D_2=\{x:\, d_G(x,S_1 \cup S_2)=2\}$.
\item No vertex in  $S_1$ is adjacent to a vertex in $S_2$, and no vertex in  $D_{11}$ is adjacent to a vertex in $D_{12}$.
\item Each vertex in $D_{11}$ is adjacent to exactly one vertex in $S_1$, each vertex of $D_{12}$ is adjacent to exactly one vertex in $S_2$,
and each vertex in $D_2$ is adjacent to exactly one vertex in $D_{11}$.
\end{enumerate}
Our main result is the following theorem which gives the characterization described above.

\begin{thm} \label{thm:girth15characterization}
If $G$ is a nontrivial graph having girth at least $15$, then  $G \in \cU$ if and only if  $G \in \cF$.
\end{thm}

The remainder of the paper is organized as follows.  In the next section we give the necessary definitions and notation used throughout the remainder of the paper.
In Section~\ref{sec:prelims}, we present a construction to prove that every graph is an induced subgraph of some graph in $\cU$.  In addition, we establish a useful
connection between $\cU$ and the class of well-covered graphs.  Section~\ref{sec:lemmas} is devoted to establishing some necessary conditions for any graph in $\cU$
that has girth at least $15$.  The structural characterization (Theorem~\ref{thm:girth15characterization}) is proved in Section~\ref{sec:theproof}, and we conclude
with some open problems in Section~\ref{sec:openproblems}.

\section{Definitions and Notation} \label{sec:defs}

Let $G$ be a finite, simple graph with vertex set $V(G)$ and edge set $E(G)$.  For $v \in V(G)$, the \emph{open neighborhood} of $v$ is the set
$N_G(v)$ of all vertices in $G$ that are adjacent to $v$.  For a subset $S$ of $V(G)$ the open neighborhood of $S$ is denoted by $N_G(S)$ and is defined
by $N_G(S)=\cup_{x\in S}N_G(x)$.  The \emph{closed neighborhood} of $v$ is the set $N_G[v]$ defined by $N_G[v]=N_G(v) \cup \{v\}$.
Whenever the graph $G$ is understood from the context, it will be removed from the subscript.  For $A \subseteq V(G)$, the subgraph of $G$ induced
by $A$ will be denoted by $G[A]$.  The \emph{independence number} of $G$ is the
maximum cardinality, $\alpha(G)$, of a set of vertices that are pairwise non-adjacent, and the  \emph{independent domination number} of $G$ is the smallest cardinality
$i(G)$ of such an independent set that is maximal in the subset inclusion relation.  A graph is  \emph{well-covered} (see Plummer~\cite{pl-1970}) if all of its maximal independent sets have the same cardinality.  That is, a graph $G$ is well-covered if $i(G)=\alpha(G)$.

A vertex $x$ of $G$ is a \emph{leaf} if $\deg(x)=1$, and
a vertex is called a \emph{support vertex} if it is adjacent to at least one leaf.  A support vertex that has more than one leaf as a neighbor is called
a \emph{strong support} vertex.  We denote the set of all leaves (support vertices) of $G$ by $L_G$ (respectively $S_G$).
The set of leaves of $G$ adjacent to a support vertex $s$ is denoted $L_G(s)$.    If $u$ and $v$ are vertices in $G$, then the distance
between $u$ and $v$ is denoted $d_G(u,v)$ and is the length of a shortest $uv$-path in $G$.  For $A \subseteq V(G)$ and a vertex $u$, we
let $d_G(u,A)$ denote $ \min\{d_G(u,v):\, v \in A\}$.  The \emph{girth} of $G$, denoted $g(G)$, is the length of a shortest cycle in $G$.  If $G$ is acyclic, then we write $g(G)=\infty$.  For a positive integer $n$, we let $[n]$ be the set of integers $\{1,2,\ldots,n\}$.

Suppose $G$ is a graph with a support vertex $s$.  Fix a leaf $x \in L_G(s)$.
Let $G'$ be the graph with $V(G')=V(G) \cup \{w\}$ and $E(G')=E(G) \cup \{sw\}$, where $w$ is a new vertex.  If $A$ is a maximal open packing of $G$, then
it follows that $A$ is a maximal open packing of $G'$.  Furthermore, if there exists $z \in A \cap L_G(s)$, then $(A-\{z\}) \cup \{w\}$ is a maximal open
packing of $G'$.  If $B$ is a maximal open packing of $G'$ and $w \notin B$, then $B$ is a maximal open packing of $G$.  On the other hand, if $w \in B$,
then $(B-\{w\})\cup \{x\}$ is a maximal open packing of $G$.  This gives the following observation.

\begin{ob} \label{ob:addingleaves}
If $G$ is a graph with at least one support vertex and $G \in \cU$, then any supergraph of $G$ that is obtained from $G$ by adding a new leaf adjacent to
any support vertex of $G$ is also in $\cU$.
If $G$ is a graph with a strong support vertex $s$ and $G \in \cU$, then $G-w \in \cU$ where $w \in L_G(s)$.
\end{ob}

\section{Preliminary Results} \label{sec:prelims}

For a given graph $G$,  the \emph{open neighborhood graph of $G$} is the graph, $\ong(G)$, whose vertex set is $V(G)$ such that distinct vertices
$u$ and $v$ are adjacent in  $\ong(G)$  if and only if $N_G(u) \cap N_G(v) \neq \emptyset$.  It is clear that a subset $A \subseteq V(G)$ is a
(maximal) open packing in $G$ if and only if $A$ is a (maximal) independent set in the graph $\ong(G)$.  Consequently,
$i(\ong(G))=\Lopack(G)\le \opack(G)=\alpha(\ong(G))$.
Thus we have the following connection between well-covered graphs and the class $\cU$.

\begin{prop} \label{prp:equivalence}
A graph $G$ is in $\cU$ if and only if $\ong(G)$ is well-covered.
\end{prop}
 A straightforward analysis shows that for $n \ge 3$ we have $\ong(C_{2n})=2 C_n$ and
$\ong(C_{2n-1})= C_{2n-1}$, while $\ong(C_3)=C_3$ and $\ong(C_4)=2P_2$.  In addition, $\ong(P_{2n})=2 P_n$ and $\ong(P_{2n+1})= P_n \cup P_{n+1}$ for every positive integer $n$, while $\ong(P_{1})= P_1$.  The next result then follows by using Proposition~\ref{prp:equivalence} and  what is known about well-covered cycles and paths.

\begin{prop}\label{prp:cycles-paths}
If $n$ is a positive integer, then
\begin{enumerate}[(i)]
\item $P_n \in \cU$ if and only if $n\in\{1,2,3,4,8\}$.
\item $C_n \in \cU$ if and only if $n \in \{3,4,5,6,7,8,10,14\}$.
\end{enumerate}
\end{prop}

The next proposition shows that, regardless of girth, there does not exist a forbidden subgraph characterization for the graphs in $\cU$.

\begin{prop} \label{prp:noforbidden}
If $H$ is any graph, then there exists a graph $G \in \cU$ such that $H$ is an induced subgraph of $G$.
\end{prop}
\proof
Suppose $H$ has order $n$ with vertex set $\{h_1,h_2,\ldots,h_n\}$.  For each $i \in [n]$, let $a_ib_ic_i$ be a path of order $3$.  A graph $G$
of order $4n$ is now constructed from the disjoint union of $H$ and the $n$ disjoint paths of order $3$ by adding the set of edges,
$\{h_ia_i:\, i\in [n]\}$.  Suppose $P$ is any maximal open packing in $G$.  It is clear that $1 \le |P \cap \{h_i,a_i,b_i,c_i\}|\le 2$, for every $i \in [n]$.
Suppose there exists $k \in [n]$  such that $|P \cap \{h_k,a_k,b_k,c_k\}|=1$, say $P \cap \{h_k,a_k,b_k,c_k\}=\{x\}$.  Each of the four possibilities for $x$
leads to the contradiction that $P$ is not maximal.  If $x\in \{h_k,b_k\}$, then $c_k$ can be added to $P$.  On the other hand,
if $x\in \{a_k,c_k\}$, then $b_k$ can be added to $P$.  This implies that $|P|=2n$ and therefore $G \in \cU$.  \qed

\section{Necessary Conditions} \label{sec:lemmas}

In this section we will derive a list of necessary conditions that are true about any graph $G \in \cU$ that has girth at least $15$.
In Section~\ref{sec:theproof} we use these conditions to prove the main characterization theorem.  Note that some of these do not require such a restriction on the girth.

\begin{lem} \label{lem:supports}
Let $G$ be a triangle-free graph.  If  $G \in \cU$, then no vertex of $G$ is adjacent to more than one support vertex.
\end{lem}
\proof
Suppose that $x$ is a vertex in a graph $G$ such that $x$ is adjacent to at least two support vertices, say $s_1$ and $s_2$.  Let $y_1$ and $y_2$ be vertices of
degree $1$ adjacent to $s_1$ and $s_2$ respectively.  Extend $\{x,s_1\}$  to a maximal open packing $P$.  Note that $P \cap (N(\{x,s_1,s_2\})-\{x,s_1\})=\emptyset$.
Since $s_1s_2 \notin E(G)$ ($G$ is triangle-free), it follows that $(P-\{x\}) \cup \{y_1,y_2\}$ is an open packing, which implies that
\[\Lopack(G)\le |P|<|(P-\{x\}) \cup \{y_1,y_2\}| \le \opack(G)\,,\]
and therefore $G \notin \cU$. \qed

\begin{lem} \label{lem:girthatleast15}
Let $G$ be a graph with girth at least $15$. If $\delta(G) \ge 2$, then $G \notin \cU$.
\end{lem}
\proof
Suppose $G$ is a graph with girth at least $15$ and $\delta(G)\ge 2$.  If also $\Delta(G) \le 2$, then $G \notin \cU$ by Proposition~\ref{prp:cycles-paths}.
Thus, we assume that $\Delta(G) \ge 3$.
Let $x$ be a vertex of degree at least $3$; suppose $N(x)=\{y_1,y_2,\ldots,y_k\}$ for some $k \ge 3$.  Consider a breadth-first search
spanning tree $T$ of $G$ rooted at $x$.  For each $i\in [7]$, let $L_i=\{u \in V(G) :\, d_G(x,u)=i \}$.  Since $G$ has girth at least $15$, we
see that $L_k$ is an independent set for $1 \le k \le 6$.  We now define a subset $A \subseteq V(G)$, beginning with $A=\emptyset$.
There are a number of situations to consider, and we will add vertices to $A$ as follows.
\begin{itemize}
\item For each $u \in L_3$ such that $d_G(y_1,u)=2$, choose vertices $u_4\in N(u) \cap L_4$ and $u_5 \in N(u_4) \cap L_5$.  Add $u_4$ and $u_5$ to $A$.

\item Select a single vertex $w \in N(y_2) \cap L_2$.  For each vertex $p \in L_4$ such that $d_G(p,w)=2$, choose vertices
$p_5\in N(p) \cap L_5$ and $p_6 \in N(p_5) \cap L_6$.  Add $p_5$ and $p_6$ to $A$.

\item For each $s \in (N(y_2)\cap L_2)-\{w\}$ and for each  $s' \in N(s) \cap L_3$, select $s_4'\in N(s')\cap L_4$ and $s_5' \in N(s_4') \cap L_5$.
Add $s_4'$ and $s_5'$ to $A$.

\item Select a single vertex $z \in N(y_3) \cap L_2$.  For each vertex $q \in L_4$ such that $d_G(q,z)=2$, choose vertices
$q_5\in N(q) \cap L_5$ and $q_6 \in N(q_5) \cap L_6$.  Add $q_5$ and $q_6$ to $A$.

\item  For each $t \in (N(y_3)\cap L_2)-\{z\}$ and for each  $t' \in N(t) \cap L_3$, select $t_4'\in N(t')\cap L_4$ and $t_5' \in N(t_4') \cap L_5$.
Add $t_4'$ and $t_5'$ to $A$.
\item For each $j\in [k]-\{1,2,3\}$ and for each vertex $r\in N(y_j)\cap L_2$ select one vertex $r_3\in N(r)\cap L_3$ and one vertex $r_4\in N(r_3)\cap L_4$.
Add $r_3$ and $r_4$ to $A$.
\end{itemize}

Note that the set $A$ constructed above is an open packing in $G$.  Let $H$ be the subgraph of $G$ induced by $S=\{x,y_1,y_2,y_3,w,z\}$.  Extend
$A$ to a maximal open packing, $B$, of the induced subgraph $G-S$ of $G$.  By the choice of the vertices placed into $A$, we see that if $g \in V(G-S)$
and $g$ is within distance $2$ of a vertex of $S$, then $g \notin B$.  On the other hand,
if a maximal open packing of $H$ is added to $B$, then the resulting set is a maximal open packing of $G$.  This implies that both
$B \cup \{x,y_2\}$ and $B \cup \{y_1,w,z\}$ are maximal open packings of $G$.  Therefore, $G \notin \cU$. \qed

\begin{lem} \label{lem:closetosupport}
Let $G$ be a graph with $\delta(G)=1$ and with girth at least $11$.  If $G \in \cU$, then every vertex of $G$ is within distance
$2$ of a support vertex.
\end{lem}
\proof
Suppose, for the sake of contradiction, that there exists $G \in \cU$ such that $\delta(G)=1$ and $g(G) \ge 11$, but for some vertex $w$ of $G$,
the distance from $w$ to the nearest support vertex is at least $3$.  Let $s$ be a support vertex of $G$ that is closest to $w$ and let $r$
be a vertex of degree $1$ adjacent to $s$.  Suppose $v$ is a vertex on a shortest $w,s$-path such that $d_G(v,s)=3$.  Let $vxzs$ be a shortest
$v,s$-path.   For each  $i\in [5]$, let $L_i=\{u \in V(G):\, d_G(v,u)=i\}$.  Since $g(G) \ge 11$, the set $L_i$ is independent for each $i \in [4]$
and no vertex in $L_5$ has more than one neighbor in $L_4$.
In addition, since no support vertex of $G$ is within distance less than $3$ of $v$, it follows that $N(u)\cap L_4\neq \emptyset$, for each $u \in L_3$.
Let $N(v)-\{x\}=\{x_1,\ldots,x_k\}$.  Let $j\in [k]$.  For each $a \in N(x_j)-\{v\}$, choose $a'\in L_3\cap N(a)$ and $a'' \in L_4 \cap N(a')$.
Let
\[ D=\bigcup_{j=1}^k \left(\bigcup\left\{a',a'':\, a \in N(x_j)-\{v\}\right\}\right)\,.\]
Since $g(G) \ge 11$, it follows that $D \cup \{z,s\}$ is an open packing of $G$ that can be extended to a maximal open packing, $P$, of $G$.
Note that $P \cap [\{x_1,\ldots,x_k\}]=\emptyset$.
However, now $(P-\{z\}) \cup \{r,v\}$ is a larger open packing of $G$, which is a contradiction.  \qed

\begin{lem} \label{lem:girthatleast7}
Let $G$ be a graph with girth at least $7$.  If there exists a path $s_1u_1vu_2s_2$ in $G$ such that $s_1$ and $s_2$ are support vertices and such
that for each $i\in [2]$, the vertex $u_i$ is the only neighbor of $s_i$ that has degree at least $2$, then $G \notin \cU$.
\end{lem}
\proof
For $i \in [2]$, let $x_i$ be a vertex of degree $1$ such that $x_i \in N(s_i)$.  In addition, define the following sets of
vertices.
\begin{itemize}
\item $A=N(u_1)-\{v,s_1\}$, $B= N(v)-\{u_1,u_2\}$, $C=N(u_2)-\{v,s_2\}$,
\item $A'=N(A)-\{u_1\}$, $B'=N(B)-\{v\}$, $C'=N(C)-\{u_2\}$.
\end{itemize}
Note that the sets $A,B,C,A',B',C'$ are pairwise disjoint since $g(G) \ge 7$.  Extend the open packing $\{x_1,v,u_2\}$ to a maximal open
packing, $P$, of $G$.  It follows that $P \cap (A \cup \{s_1,u_1\})=\emptyset$, $P \cap (B\cup C\cup B'\cup C')=\emptyset$, and $P \cap (N[s_2]-\{u_2\})=\emptyset$.
The set $Q=(P-\{v,u_2\}) \cup \{s_1,s_2,x_2\}$ is an open packing, and $|Q|>|P|$.  Therefore, $G \notin \cU$. \qed

By Lemma~\ref{lem:supports} it follows that for a triangle-free graph $G$ in $\cU$, if $G$ has vertices of degree $1$, then the subgraph of $G$ induced by
$S_G$ is a disjoint union of isolated vertices and edges.  If $s$ is such an isolated vertex in $G[S_G]$, then $s$ is called a \emph{single star support} vertex.
If $uv$ is an edge in $G[S_G]$, then $u$ and $v$ are called \emph{double star supports}.

\begin{lem} \label{lem:singlestarsupports}
Let $G$ be a connected graph such that $\delta(G)=1$ and $g(G) \ge 15$.  If $G \in \cU$ and $s$ is a single star support in $G$, then $s$ has at most
one neighbor that does not belong to $L_G$.
\end{lem}
\proof
Note that if $u$ is a single star support vertex and $v \in S_G$ such that $u \neq v$, then $d_G(u,v) \ge 3$ since by Lemma~\ref{lem:supports} no vertex of $G$
is adjacent to two support vertices.  We will prove the lemma by establishing a sequence of claims.

\begin{claim} \label{claim:1}
 If $s_1$ and $s_2$ are both single star support vertices, then $d_G(s_1,s_2) \neq 3$.
\end{claim}
\proof
Suppose for the sake of contradiction that there exists a path $s_1abs_2$ in $G$.  Let $k_i$ be a leaf adjacent to $s_i$, for $i\in [2]$.
For $i \in [4]$, let
\[L_i(s_1)=\{u:\, d_G(s_1,u)=i \text{ such that no shortest } us_1\text{-path contains }a\}-L_G\,,\]
and
\[L_i(s_2)=\{u:\, d_G(s_2,u)=i \text{ such that no shortest } us_2\text{-path contains }b\}-L_G\,.\]
For each $x \in L_1(s_1)$ and for each $x_2 \in N(x)\cap L_2(s_1)$ choose a vertex $x_3 \in N(x_2) \cap L_3(s_1)$ and a vertex $x_4 \in N(x_3) \cap L_4(s_1)$.
Similarly, for each $x \in L_1(s_2)$ and for each $x_2 \in N(x)\cap L_2(s_2)$ choose a vertex $x_3 \in N(x_2) \cap L_3(s_2)$ and a vertex $x_4 \in N(x_3) \cap L_4(s_2)$.
Let
\[P=\{x_3:\, x\in L_1(s_1)\cup L_1(s_2)\} \cup \{x_4:\, x\in L_1(s_1)\cup L_1(s_2)\} \cup \{a,b\}\,.\]
Since the girth of $G$ is more than $13$, the set $P$ is an open packing.  Extend $P$ to a maximal open packing $Q$ of $G$.  Since $(Q-\{a,b\})\cup \{s_1,k_1,s_2,k_2\}$ is also an open packing, we have reached a contradiction.  This proves Claim~\ref{claim:1}.

\begin{claim} \label{claim:2}
If $s_1$ is a single star support vertex, then there does not exist a double star support whose distance to $s_1$ is exactly $3$.
\end{claim}
\proof
Suppose the claim is false.  Let $s_1$ be a single star support and let $s_2$ and $s_3$ be adjacent double star supports such that $s_1abs_2s_3$
is a path in $G$.  For each $i \in [3]$, let $k_i \in L_G \cap N(s_i)$.  For each $i \in [4]$, let $L_i(s_1)$ be defined as in the proof of Claim~\ref{claim:1}
and let
\[L_i(s_2)=\{u:\, d_G(s_2,u)=i \text{ such that no shortest } us_2\text{-path contains }b \text{ or } s_3\}-L_G\,.\]
For each $x\in L_1(s_1) \cup L_1(s_2)$ choose $x_3$ and $x_4$ as in the proof of Claim~\ref{claim:1}.  Let
\[P=\{x_3:\, x\in L_1(s_1)\cup L_1(s_2)\} \cup \{x_4:\, x\in L_1(s_1)\cup L_1(s_2)\} \cup \{a,b,k_3\}\,.\]
By the girth assumption on $G$, it follows that $P$ is an open packing in $G$.
Extend $P$ to a maximal open packing $Q$ of $G$.  Since $(Q-\{a,b\})\cup \{s_1,k_1,k_2\}$ is also an open
packing, we have reached a contradiction.  This proves Claim~\ref{claim:2}.

\begin{claim} \label{claim:3}
No pair of single star support vertices are at distance exactly $4$.
\end{claim}
\proof
Suppose the claim is not true.   Let $s_1$ and $s_2$ be single star support vertices  such that $d_G(s_1,s_2)=4$.  Let $s_1abcs_2$ be
a path in $G$ and let $k_i \in L_G \cap N(s_i)$ for $i \in [2]$.  In a manner similar to the proofs of the above claims, for each $i \in [4]$, we let
\[L_i(s_1)=\{u:\, d_G(s_1,u)=i \text{ such that no shortest } us_1\text{-path contains }a\}-L_G\,,\]
and
\[L_i(s_2)=\{u:\, d_G(s_2,u)=i \text{ such that no shortest } us_2\text{-path contains }c\}-L_G\,.\]
Also, $P=\{x_3:\, x\in L_1(s_1)\cup L_1(s_2)\} \cup \{x_4:\, x\in L_1(s_1)\cup L_1(s_2)\} \cup \{a,b,k_2\}$.
Since $g(G) \ge 15$, the set $P$ is an open packing in $G$, and we extend it to a maximal open packing, $Q$, in $G$.
However, $(Q-\{a,b\}) \cup \{k_1,s_1,s_2\}$ is a larger open packing.  This contradiction proves the claim.

Now, let $s$ be a single star support vertex that has distinct neighbors $b$ and $c$ that both
have degree larger than $1$.  Since $s$ is a single star support vertex, it follows by definition that neither $b$ nor $c$ is
a support vertex.  Fix $a \in N(b)-\{s\}$ and $d \in N(c)-\{s\}$.  By Lemma~\ref{lem:supports}, neither $a$ nor $d$ is a support
vertex.  For each $i \in [4]$, let
\[L_i(a)=\{u:\, d_G(a,u)=i \text{ such that no shortest } ua\text{-path contains }b\}\,,\]
and
\[L_i(d)=\{u:\, d_G(d,u)=i \text{ such that no shortest } ud\text{-path contains }c\}\,.\]
By Claim~\ref{claim:1} and Claim~\ref{claim:2}, $L_1(a) \cup L_1(d)$ does not contain a support vertex, which implies
that $L_2(a) \cup L_2(d)$ does not contain any leaves.  In addition, Claim~\ref{claim:3} implies that $L_2(a) \cup L_2(d)$
does not contain any single star support vertices.  Taken together this means that each vertex in $L_2(a)$ has a neighbor in $L_3(a)$
which in turn is adjacent to a vertex in $L_4(a)$.  A similar conclusion holds for each vertex in $L_2(d)$.
For each $x \in L_1(a)$ and for each $x_2 \in N(x)\cap L_2(a)$ choose a vertex $x_3 \in N(x_2) \cap L_3(a)$ that is adjacent
to a vertex $x_4 \in L_4(a)$.  In an analogous way for
each $x \in L_1(d)$ and for each $x_2 \in N(x)\cap L_2(d)$ choose a vertex $x_3 \in N(x_2) \cap L_3(d)$ that is adjacent
to a vertex $x_4 \in L_4(d)$.  By the choice of these vertices and because the girth of $G$ is at least $15$, we see that
\[P=\{x_3:\, x\in L_1(a)\cup L_1(d)\} \cup \{x_4:\, x\in L_1(a)\cup L_1(d)\} \cup \{b,s\}\,\]
is an open packing in $G$.  We extend $P$ to a maximal open packing, $Q$, in $G$.  However, $(Q-\{s\})\cup \{a,d\}$ is a larger
open packing.  This final contradiction proves the lemma.  \qed

\begin{lem} \label{lem:distance2}
Let $G$ be a connected graph such that $\delta(G)=1$ and $g(G) \ge 15$.  If $G \in \cU$, then every vertex at distance $2$ from $S_G$ has exactly one single
star support vertex at distance exactly $2$.
\end{lem}
\proof
Suppose $G$ is a graph of girth at least $15$ such that $G \in \cU$ and $\delta(G)=1$.  Let $v$ be a vertex such that $d_G(v,S_G)=2$.  By Claim~\ref{claim:3} in
the proof of Lemma~\ref{lem:singlestarsupports}, there cannot be more than one single star support vertex at distance exactly $2$ from $v$.  Suppose there are none.
This means that there exist double star support vertices, say $s_1$ and $s_2$, and a shortest path $vas_2s_1k_1$, where $k_1\in L_G \cap N(s_1)$.  For each
$i \in [4]$, let
\[L_i(v)=\{u:\, d_G(v,u)=i \text{ such that no shortest } uv\text{-path contains }a\}\,.\]
Let $x \in L_1(v)$ and let $x_2 \in N(x) \cap L_2(v)$.  Since $x \notin S_G$, we get $\deg(x_2) \ge 2$.  By our assumption that no single star support vertex has
distance exactly $2$ from $v$, we infer that there exist $x_3 \in N(x_2) \cap L_3(v)$ and  $x_4 \in N(x_3) \cap L_4(v)$.  Let
\[P=\{x_3:\, x\in L_1(v)\} \cup \{x_4:\, x\in L_1(v)\} \cup \{s_1,s_2\}\,.\]
Since $g(G) \ge 15$, the set $P$ is an open packing.  Extend $P$ to a maximal open packing, $Q$, of $G$.  This leads to a contradiction
since $(Q-\{s_2\}) \cup \{v,k_1\}$ is an open packing of cardinality larger than $|Q|$.  \qed

\section{Proof of Theorem 1} \label{sec:theproof}

In this section we prove Theorem~\ref{thm:girth15characterization}.  For this purpose we introduce the following notation.
The set of single star support vertices in $G$  will be denoted by $S_1(G)$, and $S_2(G)$ will denote the set of double star support vertices in $G$.
By Lemma~\ref{lem:closetosupport}, every $v$ in a graph of girth at least $15$ that belongs to $\cU$ is within distance $2$ of a support vertex.  Hence,
we define $D_1(G)=\{ x:\, d_G(x,S_G)=1\}-L_G$ and $D_2(G)=\{ x:\, d_G(x,S_G)=2\}$.  For simplification when the graph is clear from the context we simply write $S_1$ instead of $S_1(G)$, and so on. By Lemma~\ref{lem:supports}, the set $D_1$ further partitions into $D_{11}$ and $D_{12}$ defined by
$D_{11}=\{ x \in D_1:\, N(x)\cap S_1 \neq \emptyset\}$ and $D_{12}=\{ x \in D_1:\, N(x)\cap S_2 \neq \emptyset\}$.

We restate Theorem~\ref{thm:girth15characterization} for convenience.

\noindent \textbf{Theorem~\ref{thm:girth15characterization}} \emph{
If $G$ is a nontrivial graph having girth at least $15$, then  $G \in \cU$ if and only if  $G \in \cF$.
}
\proof
Suppose $G$ is in $\cU$ and $G$ has girth at least $15$. It follows from Lemma~\ref{lem:girthatleast15} that $\delta(G)=1$.
Applying Lemmas~\ref{lem:supports},~\ref{lem:closetosupport},~\ref{lem:girthatleast7},~\ref{lem:singlestarsupports}, and~\ref{lem:distance2}
we see that $G \in \cF$.

Now suppose $G$ belongs to the family $\cF$ and let $L, S_1,S_2,D_{11}, D_{12},D_2$ be a weak partition of $V(G)$ that satisfies conditions
1-5 in the definition of $\cF$.  Let $P$ be any maximal open packing in $G$.  We claim that $|P|=2|S_1|+|S_2|$.  Let $S_1=\{a_1,\ldots,a_n\}$ and
let $S_2=\{b_1,\ldots,b_m,c_1,\ldots,c_m\}$ where $b_ic_i \in E(G)$ for each $i\in [m]$.  For $i \in [n]$, let $u_i$ be a leaf adjacent to $a_i$,
let $X_i=\{x \in V(G):\, d_G(x,a_i)\le 2\}$, and let $A_i=X_i \cap D_2$.  For $j \in [m]$, let $v_j \in L \cap N(b_j)$, let $w_j \in L \cap N(c_j)$, let
$Y_j=\{x \in V(G):\, d_G(x,\{b_j,c_j\})\le 1\}$, and let $B_j=Y_j \cap D_{12}$.  It is clear that $|P \cap X_i| \le 2$ and that
$|P \cap Y_j| \le 2$ for each $i\in [n]$ and each $j \in [m]$.

Suppose first that there exists $ p \in [n]$ such that $|P\cap X_p|=1$.  For the sake of reference, let $\{x\}=P\cap X_p$.  Since $P$ is
a maximal open packing, we have a contradiction.  Indeed, if $x \in L \cap X_p$ or $x \in D_{11} \cap X_p$, then $P$ can be expanded to include $a_p$.
On the other hand,  if $x=a_p$ or $x \in A_p$, then $P$ can be expanded to include $u_p$.  Now  suppose there is $q \in [m]$ such that $|P\cap Y_q|=1$,
say $\{y\}= P\cap Y_q$.  Again we arrive at a contradiction since $P$ is a maximal open packing.  In particular, if $y \in L \cap N(b_q)$, $y=c_q$,
or  $y\in D_{12} \cap N(b_q)$, then $P$ can be expanded to include $w_q$.  On the other hand if $y \in L \cap N(c_q)$, $y=b_q$, or
$y\in D_{12} \cap N(c_q)$,  then $P$ can be expanded to include $v_q$.  Therefore, $|P|=2|S_1|+|S_2|$, and it follows that $G \in \cU$.  \qed

\begin{figure} [ht!]
\begin{center}
\begin{tikzpicture}[scale=.55,style=thick] %
\def\vr{5pt}
\path (4,4) coordinate (a1); \path (6,4) coordinate (a2); \path (8,4) coordinate (a3); \path (10,4) coordinate (u1);
\path (12,4) coordinate (u2); \path (14,4) coordinate (v1); \path (15,2) coordinate (aa1); \path (15,-2) coordinate (aa2);
\path (14,-4) coordinate (vv3); \path (12,-4) coordinate (vv2); \path (10,-4) coordinate (uu5); \path (8,-4) coordinate (uu4);
\path (6,-4) coordinate (uu3); \path (4,-4) coordinate (uu2); \path (3,-2) coordinate (vvv2); \path (3,2) coordinate (uuu2);
\path (1,-2) coordinate (vvv); \path (1,2) coordinate (uuu); \path (0,2) coordinate (uuu1); \path (0,-2) coordinate (vvv1);
\path (7,-6) coordinate (uu); \path (13,-6) coordinate (vv); \path (7,-7) coordinate (uu1); \path (13,-7) coordinate (vv1);
\path (6,6) coordinate (a); \path (6,7) coordinate (b); \path (6,8) coordinate (c);
\path (11,6) coordinate (u); \path (14,6) coordinate (v); \path (11,7) coordinate (u3); \path (14,7) coordinate (v2);
\path (17,0) coordinate (aa); \path (18,0) coordinate (bb); \path (19,0) coordinate (cc);
\path (7,-6) coordinate (uu); \path (7,-7) coordinate (uu1); \path (13,-6) coordinate (vv); \path (13,-7) coordinate (vv1);

\foreach \i in {1,...,3}
{  \draw (a\i) -- (a) ; \draw (u\i) -- (u) ; \draw (vv\i) -- (vv) ;}
\foreach \i in {1,...,2}
{   \draw (v\i) -- (v) ; \draw (aa\i) -- (aa) ;  \draw (uuu\i) -- (uuu) ; \draw (vvv\i) -- (vvv) ;}
\foreach \i in {1,...,5}
{  \draw (uu\i) -- (uu) ; }
\draw (u)--(v); \draw (uu)--(vv); \draw (uuu)--(vvv); \draw (a)--(b)--(c); \draw (aa)--(bb)--(cc);

\draw (a2) -- (a3); \draw (a3) -- (u1); \draw (u1) -- (u2); \draw (u1) -- (uu5); \draw (u1) -- (uu2); \draw (u2) -- (v1);
\draw (u2) -- (uu5); \draw (u2) -- (uu4); \draw (v1) -- (aa1); \draw (aa1) -- (uuu2); \draw (aa2) -- (vv3); \draw (aa2) -- (uu5);
\draw (vv3) -- (vv2); \draw (vv2) -- (uu5); \draw (uu5) -- (uu4); \draw (vvv2) -- (uuu2);
\draw (a1) -- (uu3); \draw (a1) -- (vvv2);

\draw [fill=white] (5.85,5.85) rectangle (6.15,6.15);
\draw [fill=white] (16.85,-0.15) rectangle (17.15,0.15);
\draw [fill=black] (3.85,3.85) rectangle (4.15,4.15);  
\draw [fill=black] (5.85,3.85) rectangle (6.15,4.15);  
\draw [fill=black] (7.85,3.85) rectangle (8.15,4.15); 
\draw [fill=black] (14.85,1.85) rectangle (15.15,2.15); 
\draw [fill=black] (14.85,-2.15) rectangle (15.15,-1.85); 
\draw [fill=gray] (5.85,6.85) rectangle (6.15,7.15); 
\draw [fill=gray] (17.85,-0.15) rectangle (18.15,.15); 

\foreach \i in {1,...,2}
{  \draw (u\i)  [fill=black] circle (\vr); }
\foreach \i in {2,3}
{ \draw (vv\i)  [fill=black] circle (\vr); }
\foreach \i in {2,3,4,5}
{ \draw (uu\i)  [fill=black] circle (\vr); }
\draw (v1) [fill=black] circle (\vr); \draw (vvv2) [fill=black] circle (\vr); \draw (uuu2) [fill=black] circle (\vr);
 \draw (u) [fill=gray] circle (\vr); \draw (v) [fill=gray] circle (\vr);
\draw (vv) [fill=gray] circle (\vr); \draw (uu) [fill=gray] circle (\vr); \draw (vvv) [fill=gray] circle (\vr); \draw (uuu) [fill=gray] circle (\vr);
\draw (c) [fill=white] circle (\vr);  \draw (cc) [fill=white] circle (\vr);
\draw (u3) [fill=white] circle (\vr);\draw (v2) [fill=white] circle (\vr);\draw (uu1) [fill=white] circle (\vr);
\draw (vv1) [fill=white] circle (\vr);\draw (uuu1) [fill=white] circle (\vr);\draw (vvv1) [fill=white] circle (\vr);

\end{tikzpicture}
\end{center}
\caption{A graph in $\cF$}
\label{fig:example}
\end{figure}

We note that the proof of Theorem~\ref{thm:girth15characterization} shows that regardless of girth, $G$ belonging to $\cF$ is sufficient to guarantee
that $G\in \cU$.  See Figure~\ref{fig:example} for an illustration of a graph belonging to $\cF$.  The vertices belonging to $D_{11}$ are white squares, vertices in $D_{12}$ are solid circles, vertices in $D_2$ are solid squares, vertices in $S_1$ are gray squares, and vertices in $S_2$ are gray circles.

Graphs of girth at least $15$ that belong to the class $\cU$ can be recognized in polynomial time.  Specifically, suppose $g(G) \ge 15$.  Identify the set $L$
of vertices of degree $1$ and let $S=N(L)$.  Examine the subgraph of $G$ induced by  $S$ to determine if its components have order at most $2$.  From this
it is straightforward to determine if $G \in \cF$, that is, if there is a weak partition of $V(G)$ that satisfies the five conditions in Section~\ref{sec:intro}.

Using the above recognition algorithm applied to a tree $T$, one can of course decide whether $T \in \cU$.  It is also possible to give a different, but
equivalent, description of the trees that belong to $\cU$.  With that in mind, recall that a double star is a tree of diameter $3$ that has two support vertices.  Suppose we
have a finite collection $T_1,\ldots,T_n$ of double stars where the support vertices of $T_i$ are $u_i$ and $v_i$ for $i\in [n]$.  For each
$k \in [n]$ we label some of the leaves in $L_{T_k}(u_k)$ and in $L_{T_k}(v_k)$ as \emph{special} by applying exactly one of the following two rules to $T_k$.
\begin{itemize}
\item[(1)] We label at least one (possibly all) of the leaves in $L_{T_k}(u_k)$ as special.
\item[(2)] We label at least one, but not all, of the leaves in $L_{T_k}(u_k)$ as special.  Similarly, we label at least one, but not all, of the leaves in
$L_{T_k}(v_k)$ as special.
\end{itemize}
Suppose $A$ is the set of all vertices labeled as special.  By adding edges between some pairs of vertices in $A$ in such a way that
every vertex of $A$ is incident to at least one added edge and such that the resulting
graph is connected without any cycles, we obtain a tree $T$ that belongs to $\cU$.  To identify the partition $L, S_1,S_2,D_{11}, D_{12},D_2$ of $V(T)$
that exists for any graph in $\cF$, we
 note the following.  The set $L$ is the set of vertices of degree $1$ in $T$.  The set $S_2$ consists of the support vertices from those double stars $T_k$ that
 were treated by rule (2) above or were processed by rule (1) but for which not all the leaves adjacent to $u_k$ were labeled as special.  The set $S_1$
  consists of all $v_k$ such that $T_k$ was treated by rule (1) above and  for which all the vertices in  $L_{T_k}(u_k)$
 were labeled as special.  The sets $D_{11}$ and $D_{12}$ consist of non-leaf neighbors of vertices in $S_1$ and $S_2$, respectively.  The set $D_2$ is made up
 of all those vertices that belong to $L_{T_k}(u_k)$, where $T_k$ was treated by rule (1) such that all the vertices of $L_{T_k}(u_k)$ were labeled as special.

\begin{figure} [ht!]
\begin{center}
\begin{tikzpicture}[scale=.7,style=thick]
\def\vr{5pt}
\path (0,0) coordinate (x); \path (9,0) coordinate (y);
\path (2,2) coordinate (a1); \path (3,2) coordinate (b1);\path (4,2) coordinate (c1);
\path (5,2) coordinate (d1);\path (6,2) coordinate (e1);\path (7,2) coordinate (f1);
\path (2,1) coordinate (a2); \path (3,1) coordinate (b2);\path (4,1) coordinate (c2);
\path (5,1) coordinate (d2);\path (6,1) coordinate (e2);\path (7,1) coordinate (f2);
\path (2,-2) coordinate (ak); \path (3,-2) coordinate (bk);\path (4,-2) coordinate (ck);
\path (5,-2) coordinate (dk);\path (6,-2) coordinate (ek);\path (7,-2) coordinate (fk);
\path (4.5,.25) coordinate (el1); \path (4.5,-.5) coordinate (el2); \path (4.5,-1.25) coordinate (el3);
\draw (x)--(a1)--(b1)--(c1)--(d1)--(e1)--(f1)--(y);
\draw (x)--(a2)--(b2)--(c2)--(d2)--(e2)--(f2)--(y);
\draw (x)--(ak)--(bk)--(ck)--(dk)--(ek)--(fk)--(y);
\draw (a1) [fill=white] circle (\vr); \draw (b1) [fill=white] circle (\vr); \draw (c1) [fill=white] circle (\vr);
\draw (d1) [fill=white] circle (\vr);\draw (e1) [fill=white] circle (\vr);\draw (f1) [fill=white] circle (\vr);
\draw (a2) [fill=white] circle (\vr); \draw (b2) [fill=white] circle (\vr); \draw (c2) [fill=white] circle (\vr);
\draw (d2) [fill=white] circle (\vr);\draw (e2) [fill=white] circle (\vr);\draw (f2) [fill=white] circle (\vr);
\draw (ak) [fill=white] circle (\vr); \draw (bk) [fill=white] circle (\vr); \draw (ck) [fill=white] circle (\vr);
\draw (dk) [fill=white] circle (\vr);\draw (ek) [fill=white] circle (\vr);\draw (fk) [fill=white] circle (\vr);
\draw (x) [fill=white] circle (\vr); \draw (y) [fill=white] circle (\vr);
\draw (el1) [fill=black] circle (2pt); \draw (el2) [fill=black] circle (2pt); \draw (el3) [fill=black] circle (2pt);
\draw (-.7,0) node {$x$}; \draw (9.7,0) node {$y$};
\end{tikzpicture}
\end{center}
\caption{The graph $G_k$}
\label{fig:girth14}
\end{figure}
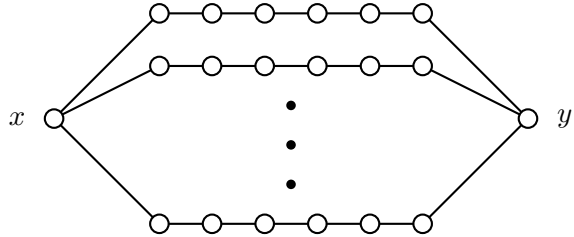

We conclude this section by noting that once the girth is lowered below $15$, more complications occur.
By Lemma~\ref{lem:girthatleast15}, every graph of girth at least $15$ that belongs to $\cU$ contains a leaf.  We now show that this girth
restriction is in some sense best possible by exhibiting an infinite family of graphs with girth $14$ and minimum degree $2$ that belong
to $\cU$.  Let $k$ be a positive integer and for each $j\in [k]$, let $a_jb_jc_jd_je_jf_j$ be a path of order $6$.   To the
disjoint union of these $k$ paths and two new vertices $x$ and $y$ we add the set $\cup_{i=1}^k\{xa_i,yf_i\}$ of edges.  The resulting graph is
denoted $G_k$.  See Figure~\ref{fig:girth14}.

We claim that $G_k \in \cU$, for each $k \ge 2$. First, note that $G_2=C_{14}\in \cU$.  Now suppose $k \ge 3$.
To verify that $G_k \in \cU$,  we appeal to Proposition~\ref{prp:equivalence} and consider $\ong(G_k)$.
Since $G_k$ is connected and bipartite, it is easy to see that $\ong(G_k)$ consists of two isomorphic components.  The component, say $H$, that contains $b_1$ has the
following structure.  The vertex set of $H$ is $\{x\} \cup \left(\cup _{i=1}^k\{b_i,d_i,f_i\}\right)$.  The set $\{f_1,f_2,\ldots,f_k\}$ induces
a clique in $H$.  The remaining edges of $H$ are those in the set $\cup_{i=1}^k \{xb_i,b_id_i,d_if_i\}$.  Let $M$ be an arbitrary maximal independent
set of $H$. If $x \in M$, then $H-N_H[x]$ is isomorphic to the corona of a clique of order $k$.  Since this corona is well-covered with independence
number $k$, it follows that $|M|=1+k$.  On the other hand if $x \not\in M$, then $M \cap \{b_1,b_2, \ldots, b_k\} \neq \emptyset$.  Without loss
of generality assume that $b_1 \in M$.  Since either $f_1$ or a neighbor of $f_1$ is in $M$ and since $\{f_1,f_2,\ldots,f_k\}$ induces a clique, it is clear that
$|M \cap \{f_1,f_2, \ldots,f_k\}|=1$.  Furthermore, $|M \cap \{b_j,d_j\}|=1$ for each $j$ such that $2 \le j \le k$.  Again we conclude that $|M|=k+1$.
Therefore, $H$ is well-covered.  By Proposition~\ref{prp:equivalence}, we infer that $G_k \in \cU$.

\section{Open Problems} \label{sec:openproblems}

We conclude with the following open problem and question.

\begin{prob}
Find a structural characterization of the class $\cU$.
\end{prob}

\begin{ques}
Is there a polynomial time algorithm to recognize the class of graphs in which all the maximal open packings have the same cardinality?
\end{ques}

\end{document}